\documentclass[12pt,a4paper]{amsart}

\usepackage {amsfonts}
\usepackage[english]{babel}
\def\g{\gamma}
\def\G{\Gamma}
\def\d{\delta}
\def\a{\alpha}
\def\b{\beta}
\def\p{\varphi}
\def\e{\varepsilon}
\def\l{\lambda}
\def\L{\Lambda}

\def\P{\Phi}

\def\R{{\mathbb R}}
\def\C{{\mathbb C}}
\def\N{{\mathbb N}}
\def\Z{{\mathbb Z}}
\def\hr{\frak r}
\def\hR{\frak R}
\def\hb{\frak b}
\def\ha{\frak a}

\def\t{\frak t}

\def\hL{\frak L}
\def\bs{~\hfill\rule{7pt}{7pt}}
\def\la{\langle}
\def\ra{\rangle}

\DeclareMathOperator{\Arg}{Arg}
\DeclareMathOperator{\supp}{supp }

\newtheorem{Th}{Theorem}
\newtheorem{Pro}{Proposition}
\newtheorem{De}{Definition}

\newtheorem{Le}{Lemma}
\begin{document}

\title{Temperate distributions with locally finite support and spectrum on Euclidean spaces}

\author{Sergii Yu.Favorov}

\address{Sergii Favorov,
\newline\hphantom{iii}  Karazin's Kharkiv National University
\newline\hphantom{iii} Svobody sq., 4, Kharkiv, Ukraine 61022}
\email{sfavorov@gmail.com}

\maketitle {\small
\begin{quote}
\noindent{\bf Abstract.}
We prove that supports of a wide class of temperate distributions with uniformly discrete support and spectrum on Euclidean spaces are finite unions of translations of full-rank lattices.
This result is a generalization of the corresponding theorem for Fourier quasicrystals, and its proof uses the technique of almost periodic distributions.
\medskip


\end{quote}
}

\medskip

The  Fourier quasicrystal may be considered as a mathematical model for  atomic arrangements
having a discrete diffraction pattern.  There are a lot of papers devoted to study properties of Fourier quasicrystals or, more generally, crystalline measures.
For example, one can mark collections of papers \cite{D}, \cite{Q}, articles \cite{F1}--\cite{F5}, \cite{K1}--\cite{M3}, and so on.

The following problem is very important in the theory of Fourier quasicrystals \cite{L1}:

{\bf Problem 1.} Let $\mu$ be a measure on $\R^d$ with discrete support $\L$ and its Fourier transform in the sense of distributions $\hat\mu$ be a measure with discrete support as well. When
is $\L$  contained in a finite union of translates of a full-rank lattice?

The answer to this question differs significantly for $d=1$ and $d>1$. In order to present the corresponding results, we need to introduce suitable notations and definitions.

         \section{Preliminaries}\label{S1}

Denote by $S(\R^d)$ the Schwartz space of test functions $\p\in C^\infty(\R^d)$ with the finite norms
 $$
  N_n(\p)=\sup_{\R^d}\max_{\|k\|\le n}\,(1+|x|)^n |D^k\p(x)|,\quad n=0,1,2,\dots,
 $$
where
$$
|x|=(x_1^2+\dots+x_d^2)^{1/2},\   k=(k_1,\dots,k_d)\in(\N\cup\{0\})^d,\ \|k\|=k_1+\dots+k_d,\  D^k=\partial^{k_1}_{x_1}\dots\partial^{k_d}_{x_d}.
 $$
 These norms generate the topology on $S(\R^d)$.  Elements of the space $S^*(\R^d)$ of continuous linear functionals on $S(\R^d)$ are called {\it temperate distributions}.
The Fourier transform of a temperate distribution $f$ is defined by the equality
\begin{equation}\label{a0}
\hat f(\p)=f(\hat\p)\quad\mbox{for all}\quad\p\in S(\R^d),
\end{equation}
where
$$
   \hat\p(y)=\int_{\R^d}\p(x)\exp\{-2\pi i\la x,y\ra\}dx
 $$
is the Fourier transform of the function $\p$. Also,
$$
   \check\p(y)=\int_{\R^d}\p(x)\exp\{2\pi i\la x,y\ra\}dx
 $$
 means the inverse Fourier transform. Note that the Fourier transform is an isomorphism of $S(\R^d)$ on $S(\R^d)$ and, respectively, $S^*(\R^d)$ on $S^*(\R^d)$.

  We will say that a set $A\subset\R^d$ is {\it locally finite} if the intersection of $A$ with any ball is finite,  $A$ is {\it relatively dense} if there is $R<\infty$
  such that $A$ intersects with each ball of radius $R$, and $A$ is {\it uniformly discrete}, if $A$ is locally finite and has a strictly positive separating constant
 $$
 \eta(A):=\inf\{|x-x'|:\,x,\,x'\in A,\,x\neq x'\}.
  $$
  A locally finite set $A$ is of {\it bounded density} if
  $$
        \sup_{x\in\R^d}\# A\cap B(x,1)<\infty.
  $$
  As usual, $\# E$ means a number of elements of the finite set $E$, and $B(x,r)$ means the ball with center in $x$ and radius $r$.

  An element $f\in S^*(\R^d)$ is called {\it a crystalline measure} if $f$ and $\hat f$ are  complex-valued measures on $\R^d$ with locally finite supports.
  The support of $\hat f$ for a  distribution $f\in S^*(\R^d)$ is called {\it spectrum} of $f$.

  Denote by $|\mu|(A)$ the variation of a complex-valued measure $\mu$ on $A$. If both measures $|\mu|$ and $|\hat\mu|$ have  locally finite supports and belong to $S^*(\R^d)$,
we say that $\mu$ is a {\it Fourier quasicrystal}. A measure $\mu=\sum_{\l\in\L}a_\l\d_\l$ with $a_\l\in\C$ and countable $\L$ is called {\it purely point}.

  Furthermore, it is well known that every distribution with locally finite support $\L$ has the form
 \begin{equation}\label{r0}
  f=\sum_{\l\in\L}\sum_k p_k(\l)D^k\d_\l,\quad k\in(\N\cup\{0\})^d,
  \end{equation}
 where the internal sum is finite for every $\l\in\L$, and $\d_\l$ is the unit mass at the point $\l$.

\begin{Pro}\label{P1}
i) If a distribution $f\in S^*(\R^d)$ has locally finite support, then
  \begin{equation}\label{r1}
f=\sum_{\l\in\L}\sum_{\|k\|\le K}p_k(\l)D^k\d_\l,\quad k\in(\N\cup\{0\})^d,
\end{equation}
where $K<\infty$ does not depend on $\l$ (\cite{F1}, Proposition 1),

ii) if a distribution $f\in S^*(\R^d)$ has uniformly discrete support, then $K$ does not depend on $\l$  and for some $M<\infty$ and all $k$
\begin{equation}\label{r2}
p_k(\l)=O(|\l|^M)\quad\mbox{as}\quad \l\to\infty,
\end{equation}
(\cite{P}, Proposition 1),

iii) if a distribution $f\in S^*(\R^d)$  of form \eqref{r1}  with locally finite support $\L$ and locally finite spectrum $\G$  has the property
\begin{equation}\label{r3}
\sum_{|\l|\le r}\sum_{\|k\|=0}^K|p_k(\l)|=O(r^{d+M})\phantom{XXXXXX}(r\to\infty)
\end{equation}
with $M\ge0$, then
\begin{equation}\label{m}
\hat f=\sum_{\g\in\G}\sum_{\|m\|\le M} q_m(\g)D^m\d_\g,\quad m\in(\N\cup\{0\})^d,
\end{equation}
If, in addition, $M$ is an integer, then for $\|m\|=M$
\begin{equation}\label{r}
q_m(\g)=O(|\g|^K)\quad\mbox{as}\quad \g\to\infty.
\end{equation}
 with the same $K$ as in \eqref{r1}. If $M$ is a positive real number, and $\G$ is uniformly discrete, then \eqref{r} is valid for all $m$ (\cite{F1}, Proposition 3),

In particular, if $f$ is a measure with uniformly bounded masses, support of bounded density, and locally finite spectrum, then $K=M=0$ and $\hat f$ is a measure with uniformly bounded masses as well.
\end{Pro}
{\bf Remark 1}. Analysis of the proof of the Proposition shows that if  we replace $O(r^{d+M})$ by $o(r^{d+M})$ in \eqref{r3},
then we may replace the inequality $\|m\|\le M$  by $\|m\|<M$ in \eqref{m}.
\medskip

{\it A lattice} is a discrete (locally finite) subgroup of $\R^d$. If $A$ is a lattice or a coset of some lattice in $\R^d$, then its rank $r(A)$ is the dimension of the smallest
translated subspace of $\R^d$ that contains $A$.  Every lattice $L$ of rank $k$ has the form $T\Z^k$, where $T:\,\Z^k\to\R^d$ is a linear operator of rank $k$. For $k=d$ we get a  full-rank lattice.
In this case the lattice
$$
L^*=\{y\in\R^d: <\l,y>\in\Z\quad \forall \l\in L\}
$$
 is called the conjugate  lattice.  It follows from Poisson's formula
 $$
 \sum_{n\in\Z^d}f(n)=\sum_{n\in\Z^d}\hat f(n),\qquad f\in S(\R^d),
 $$
   that for a full-rank lattice $L=T\Z^d$ we have
 \begin{equation}\label{a1}
      \widehat{\sum_{\l\in L} \d_\l}=|\det T|^{-1}\sum_{\l\in L^*} \d_\l.
 \end{equation}

     \section{Introduction, one-dimensional case}\label{S2}

We begin with the following result of N.Lev and A.Olevskii \cite{LO1}:
\begin{Th}\label{T1}
Let $\mu$ be a crystalline measure on $\R$ with uniformly discrete support and spectrum. Then $\supp\mu$ is a subset of a finite union of translates of a single lattice $L\subset\R$
(i.e.,of arithmetic progressions with the same difference). Moreover, $\mu$ has the form
$$
   \mu=\sum_{j=1}^N \sum_{\l\in L+\tau_j}P_j(\l)\d_\l
$$
for certain reals $\tau_j$ and trigonometric polynomials $P_j(\l)$ ($1\le j\le N$).
\end{Th}
Theorem \ref{T1} remains valid under a weaker assumption that $\supp\mu$ is a relatively dense set of bounded density (not assumed to be uniformly discrete) \cite{LO2}. However there exist
 examples of crystalline measures on $\R$, whose supports  are not contained in any finite union of translates of a lattice (see for example \cite{LO3}, \cite{M2}, \cite{KS}).

In article \cite{LR} Theorem \ref{T1} was extended to temperate distributions on the line:
\begin{Th}\label{T2}
Let $f$ be a temperate distribution on $\R$ such that $\L=\supp f$ and
$\G=\supp\hat f$ are uniformly discrete sets. Then $f$ can be represented in the form
$$
f=\sum_{\tau,\omega, l, p}c(\tau,\omega, l, p)\sum_{\l\in L}\l^le^{2\pi i\l\omega}\d_{\l+\tau}^{(p)},
$$
where $L$ is a lattice on $\R$, $(\tau,\omega, l, p)$ goes through a finite set of quadruples such that $\tau,\omega$ are
real numbers, $l, p$ are nonnegative integers, and $c(\tau,\omega, l, p)$ are complex numbers.

The result is still correct if the support $\L$ has bounded density and the coefficients $p_k(\l)$ in \eqref{r0} satisfy \eqref{r2}, while the spectrum $\G$ is
uniformly discrete.
\end{Th}

\section{Introduction, multivariate case}\label{S3}

Various analogs of Theorem \ref{T2} for temperate distributions on $\R^d$ were obtained in \cite{P} under conditions of locally finite set of differences $\L-\L$ and $\G-\G$ and in \cite{F1}
under conditions of locally finite  $\L-\L$, not too fast approaching points from $\G$, and uniformly separated from zero and infinity coefficients $p_k(\l)$.
\smallskip

In \cite{LO1}, \cite{LO2}  N.Lev and A.Olevskii  proved the following result:
\begin{Th}\label{T3}
Let $\mu$ be a {\bf positive}  measure on $\R^d$ such that $\L=\supp\mu$ and
$\G=\supp\hat\mu$ are uniformly discrete sets. Then $\L$ is contained in a finite union of translates of a lattice of rank $d$.  The same is valid if $\L$ is locally finite not assumed to be uniformly discrete.
\end{Th}
It was proved in \cite{F2} that there is a signed measure on $\R^2$ such that its support and  spectrum
are both uniformly discrete and simultaneously are unions of two incommensurable full-rank lattices. Therefore neither support, nor spectrum can be finite unions of translations of a {\bf single} lattice.

But there are results of a different type
 for  measures in $\R^d$ with discrete support, in which a {\bf finite number of  lattices} are already involved (\cite{M3},  \cite{KL}, \cite{K1}, \cite{Co}).
 In fact, the following result was proved:
\begin{Th}\label{T4}
 Let  $\mu$ be a measure on $\R^d$ with uniformly discrete support $\L$. If complex masses $\mu(\{\l\})$ at  points $\l\in\L$  take  values only from a finite set $F\subset\C$,
 and the measure $\hat\mu$ is purely point and satisfies the condition
 \begin{equation}\label{mu}
|\hat\mu|(B(0,r))=O(r^d) \quad (r\to\infty),
\end{equation}
 then $\L$  is a finite union of translations of several, possibly incommensurable, full-rank lattices.
\end{Th}
In the paper \cite{LR}, the following problem was posed:
 \smallskip

{\bf Problem 2.} Let $\mu$ be a measure on $\R^d$ with discrete support $\L$ and its Fourier transform in the sense of distributions $\hat\mu$ is a measure with discrete support as well.
Does it follow that $\L$ can be covered by a finite union of translates of several lattices?
\smallskip

Remark that the periodic structure of a crystalline measure in the multidimensional case  follows under extra assumptions on  masses (\cite{F4},\cite{F5}):
\begin{Th}\label{T5}
 Let $\mu$ be a measure on $\R^d$ with uniformly discrete support $\L$ such that $\inf_{\l\in\L}|\mu(\{\l\})|>0$ and the measure $\hat\mu$  is purely point and satisfy \eqref{mu}.
 Then  $\L$  is a finite union of translates of several disjoint full-rank lattices.

  Moreover,  there exist an integer $J$, lattices $L_1,\dots,L_J$ in $\R^d$ of rank $d$ (some of them may coincide), points $\l_1,\dots,\l_J\in\L$,
  a bounded set $\{\a_s^j:\,s\in\N,\,1\le j\le J\}\subset\R^d$, and functions
  $$
  F_j(y)=\sum_s b_s^j e^{2\pi i\la y,\a_s^j\ra}\quad\mbox{with}\quad \sum_s|b_s^j|<\infty,\quad j=1,\dots,J,
  $$
   such that
 $$
 \mu=\sum_{j=1}^J F_j(\l)\sum_{\l\in L_j+\l_j}\d_\l.
 $$
 \end{Th}
 The proofs of theorems \ref{T4} and \ref{T5} are based on Cohen's Idempotent Theorem (see, e.g., \cite{Ru};
 the first using of the Idempotent Theorem to study structure of tilings was in \cite {LM}):
\begin{Th}\label{T6}
Let $G$ be a locally compact Abelian group and $\hat G$ its dual group. If $\nu$ is a finite Borel measure on $G$ and is such that
its Fourier transform $\hat\nu(\l)$ takes only  values $0$ and $1$, then the set $\{\l:\,\hat\nu(\l)=1\}$ belongs to the coset ring of $\hat G$.
\end{Th}
Recall that the coset ring of an Abelian topological group $G$ is the smallest collection  of subsets of $G$ that is closed under complement, finite unions, and finite intersections,
  which contains all cosets of all open subgroups of $G$. In particular, if $G=\R^d_{discr}$, i.e., the space $\R^d$ with respect to the discrete topology, then the coset ring contains all cosets of all subgroups of $G$.
 \begin{Th}\label{T7}[M.Kolountzakis \cite{K1}, Theorem 3]
  All  elements of the coset ring for $\R^d_{discr}$, which are discrete in the Euclidean topology of $\R^d$, are finite unions of sets of the type $A\setminus(\cup_{j=1}^J B_j)$,
where  $A,\ B_j$ are discrete cosets in the Euclidean topology.
\end{Th}

\section{Results}\label{S4}

Here we get the following extension of Theorems \ref{T3} and \ref{T5} to temperate distributions:
\begin{Th}\label{T8}
Suppose that a temperate distribution $f$ has   locally finite support $\L$ of  bounded density, a uniformly discrete spectrum $\G$, and all coefficients $p_k(\l)$ in \eqref{r0} are
uniformly bounded and non-negative. Then $\L$ is contained in a finite union of translates of several full-rank lattices.

The non-negativity of  coefficients $p_k(\l)$  can be replaced by the following: for each $k\in(\N\cup\{0\})^d$ there are $\xi_k\in\R$ such that $\Arg\,p_k(\l)=\xi_k$
for all $p_k(\l)\neq0,\,\l\in\L$.
\end{Th}
\begin{Th}\label{T9}
Suppose that a temperate distribution $f$ of form \eqref{r0} has both uniformly discrete support $\L$ and spectrum $\G$.
If there are constants $c,\,C$ such that for  coefficients in \eqref{r0}  the inequalities  hold
\begin{equation}\label{c3}
 0<c\le\sum_k |p_k(\l)|\le C<\infty,\qquad \forall\ \l\in\L,
 \end{equation}
then $\L$ is  a  finite union of translations of several full-rank lattices.
More accurately, there exist integers $J,\,K$, lattices $L_1,\dots,L_J$ in $\R^d$ of rank $d$ (some of them may coincide), points $\l_1,\dots,\l_J\in\L$, and
a bounded set of exponents $\b _s=\b_s(j,k)\in\R^d$  such that
$$
f=\sum_{j=1}^J\sum_{\|k\|\le K}\sum_{\l\in L_j+\l_j} F_{k,j}(\l)D^k\d_\l,
$$
where
$$
F_{k,j}(x)=\sum\limits_s b_s(k,j) e^{2\pi i\la x,\b_s\ra}\quad\mbox{with}\quad\sum\limits_s|b_s(k,j)|<\infty.
$$
\end{Th}
{\bf Remark 2}.  Proposition \ref{P1} i) implies that $f$ has form \eqref{r1} with $K<\infty$, therefore we may replace the sum $\sum_k|p_k(\l)|$ by $\sup_k|p_k(\l)|$ in \eqref{c3}.
Also, we can replace the right-hand bound in \eqref{c3} by
$$
p_k(\l)=o(|\l|)\quad (\l\to\infty)\qquad\forall k,
$$
because in this case we get again $p_k(\l)=O(1)$ (see Remark 3 in Section \ref{S6}).

Nevertheless,  some upper and lower bounds are required.

Let $A=\{n\in\Z^2:\,n_1\neq0\}$. It follows from \eqref{a0} and \eqref{a1} that the Fourier transform of the measure with unbounded coefficients
 $$
   \sum_{n\in A} n_1\d_n=\sum_{x\in\Z^2}x_1\d_x
  $$
  is the distribution
$$
   \frac{i}{2\pi}\sum_{\l\in\Z^2}\frac{\partial}{\partial\l_1}\d_\l.
$$
Also, the Fourier transform of the measure whose masses do not separate from zero
$$
   \sum_{n\in A} \sin\left(\frac{2\pi n_1}{\sqrt{5}}\right)\d_n=\frac{1}{2i}\sum_{x\in\Z^2}\left(e^{2\pi ix_1/\sqrt{5}}-e^{-2\pi ix_1/\sqrt{5}}\right)\d_x
 $$
 is equal to
 $$
   \frac{1}{2i}\sum_{y\in\Z^2}\d_{y+\a}-\frac{1}{2i}\sum_{y\in\Z^2}\d_{y-\a},\qquad \a=(1/\sqrt{5},0),
$$
hence in both cases the Fourier transforms have uniformly discrete supports. But $A$  is not  a finite union of cosets of several full-rank lattices, because in the opposite case the set $\Z\setminus\{0\}$,
which is the projection of $A$ on $\R_{x_1}$, would be  a finite union of arithmetical progressions, that is impossible.
\smallskip

The above theorems  are based on the following  result for  temperate distributions with uniformly discrete spectrum, which is also of independent interest:
\begin{Th}\label{T10}
 Suppose $f\in S^*(\R^d)$  of  form \eqref{r0} has  locally finite support $\L$ of bounded density and uniformly discrete spectrum $\G$. Set for  $k\in(\N\cup\{0\})^d$
 \begin{equation}\label{c2}
\mu_k=\sum_{\l\in\L}p_k(\l)\d_\l.
\end{equation}

a) If $p_k(\l)$ are uniformly bounded in $k\in(\N\cup\{0\})^d$ and $\l\in\L$, then  $\hat\mu_k\in S^*(\R^d)$ are measures of the form
$$
 \hat\mu_k=\sum_{\g\in\G_k}q_k(\g)\d_\g
$$
 with uniformly bounded masses $q_k(\g)$  and uniformly discrete supports $\G_k$ such that
\begin{equation}\label{et}
\eta(\G_k)\ge\eta(\G)\qquad \forall k.
\end{equation}

b) if both $\L$ and $\G$ are uniformly discrete, then $\hat\mu_k\in S^*(\R^d)$ are distributions of the form
\begin{equation}\label{m2}
\hat\mu_k=\sum_{\|m\|\le M}\sum_{\g\in\G_k}q_{k,m}(\g)D^m\d_\g
\end{equation}
with uniformly bounded coefficients $q_{k,m}(\g)$ and uniformly discrete supports $\G_k$, which satisfy \eqref{et}
\end{Th}
Since $\mu_k$ is the inverse Fourier transform of $\hat\mu_k$, we see that $\mu_k\in S^*(\R^d)$. Also, by Proposition \ref{P1} i), $\mu_k\equiv0$ for $\|k\|$ large enough. Hence,
\begin{equation}\label{f2}
   f=\sum_k D^k\mu_k,\phantom{XXXXX} \hat f(y)=\sum_k (2\pi i)^{\|k\|}y^k\hat\mu_k,
 \end{equation}
where, as usually, $y^k=y_1^{k_1}\cdot\dots\cdot y_d^{k_d}$.

\section{Almost periodic functions and distributions}\label{S5}

We recall here  definitions and properties of almost periodic functions and distributions that will be used in what follows. A more complete exposition of these issues is available in \cite{C},
\cite{A}, \cite{M1}, \cite{M2}, \cite{R}.

\begin{De}\label{D1} A continuous function $g$ on $\R^d$ is  almost periodic if for any  $\e>0$ the set of its $\e$-almost periods
 $$
  \{\tau\in\R^d:\,\sup_{t\in\R^d}|g(t+\tau)-g(t)|<\e\}
 $$
is a relatively dense set in $\R^d$.
\end{De}
An equivalent definition follows:
\begin{De}\label{D2} A continuous function $g$ on $\R^d$ is  almost periodic if for any sequence $\{t_n\}\subset\R^d$ there is a subsequence  $\{t'_n\}$ such that
 the sequence of functions $g(t+t'_n)$ converge uniformly in $t\in\R^d$.
\end{De}
Using an appropriate definition, one can prove various properties of almost periodic functions.
\begin{itemize}
\item almost periodic functions are  bounded and uniformly continuous on $\R^d$,
\item the class of almost periodic functions is closed with respect to taking absolute values and linear combinations of a finite family of functions,
 \item a limit of a uniformly convergent sequence  of almost periodic functions is also almost periodic,
 \item any finite family of almost periodic functions has a relatively dense set of common $\e$-almost periods,
\item  for any almost periodic function $g(x)$ on $\R^d$ the function $h(t)=g(tx)$ is almost periodic in $t\in\R$ for any fixed $x\in\R^d$; in particular,
$g(t_1,\dots,t_d)$ is almost periodic in each variable $t_j\in\R,\ j=1,\dots,d$, if the other variables are held fixed.
 \end{itemize}
Also, we will use the following definition
\begin{De} A distribution $g$ is  almost periodic if the function $(g(y),\p(t-y))$
is almost periodic in $t\in\R^d$ for each $C^\infty$-function $\p$ on $\R^d$ with compact support.
A measure $\mu$ is almost periodic if it is an almost periodic distribution.
 \end{De}
Clearly, every linear combination of almost periodic distributions is almost periodic, and each almost periodic distribution  has a relatively dense support.

 Note that the usual definition of almost periodicity for measures (instead of $\p\in C^\infty$, we consider  continuous $\p$ with compact support)
 differs from the one given above. However these definitions coincide for nonnegative measures or measures with uniformly discrete support (see  \cite{A}, \cite{M1}, \cite{F2}).
\begin{Pro}\label{P2}
Let $\mu=\sum_{\l\in\L} p(\l)\d_\l,\  p(\l)\in\C$, be an almost periodic measure with uniformly discrete support $\L\subset\R^d$, then

i) the masses $p(\l)$ are uniformly bounded,

ii)  the measure $|\mu|=\sum_{\l\in\L}|p(\l)|\d_\l$ is almost periodic,

iii) if $\inf_{\l\in\L}|p(\l)|\ge\a>0$, then the measure $\d_\L:=\sum_{\l\in\L}\d_\l$ is almost periodic as well,

iv) for any $\e>0$ and $\l\in\R^d$ there is a relatively dense set $T$ such that $|\mu(\l')-\mu(\l)|<\e$ for all $\l'\in T$,

v) if $\mu_n,\,n=1,\dots,N$ be almost periodic measures  with uniformly discrete supports, then for any vector $\l\in\R^d$ and numbers $\e>0,\ \rho>0$ there is a relatively dense set $T$ such that
for every $\l'\in T$ and some $\l^{(n)}\in B(\l',\rho)$
\begin{equation}\label{m1}
|\mu_n(\l^{(n)})-\mu_n(\l)|<\e,\quad n=1,\dots,N.
\end{equation}
 \end{Pro}

  The condition on $\L$ to be uniformly discrete is essential (see\cite{FK}).
\smallskip

 {\bf Proof of Proposition \ref{P2}}. i) Let $\p$ be $C^\infty$ function with support in $B(0,\eta(\L)/2)$ such that $\p(0)=1$. The function $(\mu(y),\p(t-y))$ is almost periodic,
 hence it is uniformly bounded, and the numbers $p(\l)=(\mu(y),\p(\l-y))$ are uniformly bounded as well.

 ii) Every $C^\infty$-function  with compact support is a finite linear combination of shifts of nonnegative $C^\infty$-functions with supports in balls of radius $\eta/2$
 \footnote{Using Theorem on Partition of Unity, it suffices to prove this for the case of a real-valued  $\p\in C^\infty$ with support in the ball $B(x,\eta(\L)/4)$.
 Then we have $\p=\psi-(\psi-\p)$, where $\psi$ is an arbitrary $C^\infty$ nonnegative function with support in $B(x,\eta(\L)/2)$ such that $\p\le \psi$.}.
  Consequently we can check the almost periodicity  using only  such functions.

 Let $\p$ be one of them and $\supp\p\subset B(x,\eta/2)$. For each $t,\tau\in\R^d$ we have
 $$
       |(|\mu|(y),\p(t-y))-(|\mu|(y),\p(t+\tau-y))|=\left|\sum_{\l\in\L}\p(t-\l)|p(\l)|-\sum_{\l'\in\L}\p(t+\tau-\l')|p(\l')|\right|.
 $$
 Evidently, here each sum consists of at most one nonzero term. Therefore, for some $\l, \l'\in\L$ that depends on $t,\tau$ we get
 $$
       |(|\mu|(y),\p(t-y))-(|\mu|(y),\p(t+\tau-y))|=|\p(t-\l)|p(\l)|-\p(t+\tau-\l')|p(\l')||
 $$
  $$
   \le  |\p(t-\l)p(\l)-\p(t+\tau-\l')p(\l')|=|(\mu(y),\p(t-y))-(\mu(y),\p(t+\tau-y))|.
  $$
 Hence every $\e$-period $\tau$ of the function $(\mu(y),\p(t-y))$ is an $\e$-period of the function $(|\mu|(y),\p(t-y))$.

 iii) Pick $\e>0$ and $\b<\min\{\a,\eta(\L)/2\}$ such that $|\p(x)-\p(x')|<\e$ for $|x-x'|<\b$. Let $\psi$ be $C^\infty$-function such that
$\supp\psi(x)\subset B(0,\b)$ and $\psi(0)=1$.  Let  $\tau$ be some $\b$-almost period of the function $(\mu(y),\psi(t-y))$, and let $t\in B(\l,\b)$ for $\l\in\L$,
Since $\b<\eta(\L)/2$, we get for at most one $\l'\in\L$ that depends on $t+\tau$
   $$
      \b>\left|\sum_{\l\in\L}\psi(t-\l)p(\l)-\sum_{\l'\in\L}\psi(t+\tau-\l')p(\l')\right|= |p(\l)\psi(t-\l)-p(\l')\psi(t+\tau-\l')|.
   $$
For $t=\l$ we get
    $$
    |p(\l)-p(\l')\psi(\l+\tau-\l')|<\b.
    $$
 The inequality $\b<\a\le|p(\l)|$ implies that $\psi(\l+\tau-\l')\neq0$. Hence, $|\l+\tau-\l'|<\b$. Therefore, $|\p(t-\l)-\p(t+\tau-\l')|<\e$ and
  $$
       \left|\sum_{\l\in\L}\p(t-\l)-\sum_{\l'\in\L}\p(t+\tau-\l')\right|=|\p(t-\l)-\p(t+\tau-\l')|<\e.
  $$
This inequality is valid for all $\b$-almost periods $\tau$ from a relatively dense set. Thus the function $(\d_\L(y),\p(t-y))$ is almost periodic.

 iv) It is enough to consider the case $\e\le|\mu(\l)|$. Let $\rho<\eta(\supp\mu)/2$, and let $\p$ be $C^\infty$ function on $\R^d$ such that
 $$
   0\le\p(t)\le1,\quad \p(0)=1,\quad \p(t)=\p(-t),\quad \supp\p\subset B(0,\rho).
 $$
We have $\mu(\l)=(\p(\l-x),\mu(x))$.  If $\tau$ be an $\e/3$-almost period of the almost periodic function $(\p(t-x),\mu(x))$, we get
$$
|(\p(\l+\tau-x),\mu(x))-\mu(\l)|<\e/3.
$$
Therefore, there is a single point $\l'\in B(\l+\tau,\rho)$ such that
$$
|\mu(\l')\p(\l+\tau-\l')-\mu(\l)|<\e/3.
$$
We get $\l'\in\supp\mu$ and $|\mu(\l')|>2\e/3$. Since $-\tau$ is an $\e/3$-almost period as well, we see that for some $\l''\in B(\l'-\tau,\rho)$
$$
|\mu(\l'')\p(\l'-\tau-\l'')-\mu(\l')|<\e/3.
$$
We get $\l''\in\supp\mu$ and  $|\l''-\l|<2\rho$, hence, $\l''=\l$. Then we have for $t=\l+\tau-\l'$
$$
 |\mu(\l)-\mu(\l')|\le|\mu(\l)-\p(t)\mu(\l')|+|\mu(\l')-\p(t)\mu(\l')|\le|\mu(\l)-\p(t)\mu(\l')|+(1-\p^2(t))|\mu(\l')|
$$
$$
\le|\mu(\l)-\p(t)\mu(\l')|+|\mu(\l')-\p(-t)\mu(\l)|+|\p(t)\mu(\l)-\p^2(t)\mu(\l')|<\e.
$$
Recall that almost periods form a relatively dense set.

 v) The proof is very close to the previous one. We take $\rho<\min_n \eta(\mu_n)/2$,  the same function $\p$, and $\l'=\l+\tau$,
 where $\tau$ is a common $\e/3$-almost period of the functions $(\p(\l-x),\mu_n(x))$. Then, as $\l^{(n)}$ we take points satisfying  the conditions
$$
|\mu_n(\l^{(n)})\p(\l'-\l^{(n)})-\mu_n(\l)|<\e/3.
$$
Further, arguing as above, we obtain \eqref{m1}. \bs
\medskip

Typical examples of almost periodic functions on $\R^d$ are  sums of the form
\begin{equation}\label{W}
   f(t)=\sum_n a_n e^{2\pi i\la t,s_n\ra},\quad a_n\in\C,\quad s_n\in\R^d,\quad \sum_n|a_n|<\infty.
\end{equation}
It is not hard to prove that $\hat f=\sum_n a_n\d_{s_n}$.

Denote by $W$ the class of functions admitting  representation \eqref{W}. It is easy to check that a product of two functions from $W$ belongs to $W$ as well.
Also, in Section \ref{S7} we will use the following local version of the classical Wiener-Levi Theorem:
\begin{Pro}[\cite{F4}]\label{P3}
 Let  $K\subset\C$ be an arbitrary compact set,  $h(z)$ be a holomorphic function on a neighborhood of $K$, and $f\in W$. Then there is a function $g\in W$ such that if $f(x)\in K$ then $h(f(x))=g(x)$.
 \end{Pro}
Next, the Bohr compactification $\hR$ of $\R^d$ is a compact group with the dual $\R^d_{discr}$, and  $\R^d$ is a dense subset of $\hR$ with respect to the topology on $\hR$. Moreover,
restrictions to $\R^d$ of continuous functions on $\hR$ are almost periodic functions on $\R^d$ (see, e.g.,\cite{Ru}).

\section{Representation of distributions with discrete support and spectrum}\label{S6}

First we prove the following lemma:
\begin{Le}\label{L1}
Let $g_0(y), g_1(y),\dots, g_N(y)$ be temperate distributions on $\R^d$ such that for every $C^\infty$-function $\p(y)$ with compact support the functions
$$
(g_n(y),\p(t-y)),\quad t=(t_1,t_2,\dots,t_d),\quad n=0,\dots,N,
$$
are almost periodic in $t_1\in\R$ for every fixed $(t_2,\dots,t_d)\in\R^{d-1}$. If the distribution
$$
F=\sum_{n=0}^N y_1^n g_n(y),\quad y=(y_1,\dots,y_d),
$$
 has  a uniformly discrete support $T$, then the set $\cup_n\supp g_n$ is uniformly discrete and $\eta(\cup_n\supp g_n)\ge\eta(T)$.
\end{Le}

{\bf Proof}. Note that for every $C^\infty$-function $\p(y)$ with compact support we have
\begin{equation}\label{s1}
   (F(y),\p(t-y))=\sum_{n=0}^N(g_n(y),y_1^n\p(t-y))=\sum_{m=0}^N t_1^m \sum_{n=m}^N \P_{n,m}(t),
\end{equation}
where
$$
\P_{n,m}(t)=\binom{n}{m}(g_n(y),(y_1-t_1)^{n-m}\p(t-y)),\quad n\ge m.
$$
We will prove that for every $n,m\le N$ and every $y\in\supp g_n,\,\tilde y\in\supp g_m$ we get $|y-\tilde y|\ge\eta(T)$.

Consider the case $n=m=N$. Assume the converse. Then there are points $y',y''\in\supp g_N$ such that $0<|y'-y''|<\eta(T)$.
Pick
$$
\a<(1/2)\min\{|y'-y''|,\,(\eta(T)-|y'-y''|)\}
$$
 and $C^\infty$-function $\p$ with $\supp\p\subset B(0,\a)$ such that
 $$
\P_{N,N}(y')=(g_N(y),\p(y'-y))\neq0,\qquad \P_{N,N}(y'')=(g_N(y),\p(y''-y))\neq0.
 $$
 Take $\e>0$ such that $|\P_{N,N}(y')|>2\e$, $|\P_{N,N}(y'')|>2\e$. If $\tau$ is a common $\e$-almost period of the functions
 $$
    \P_{N,N}(y'+se_1),\quad \P_{N,N}(y''+se_1),\quad\mbox{where}\quad e_1=(1,0,0\dots 0),\quad s\in\R,
 $$
 then $|\P_{N,N}(y'+\tau e_1)|>\e$, $|\P_{N,N}(y''+\tau e_1)|>\e$. We have
$$
(y'_1+\tau)^{-N}(F(y),\p(y'+\tau e_1-y))=\P_{N,N}(y'+\tau e_1)+\sum_{m=0}^{N-1} (y'_1+\tau)^{m-N} \sum_{n=m}^{N-1}\P_{n,m}(y'+\tau e_1).
$$
Note that all  functions $\P_{n,m}(y'+se_1),\,0\le n,m\le N,$ are almost periodic in $s$, hence they are uniformly bounded on $\R$. Therefore if $\tau$ is large enough then
$$
(y'_1+\tau)^{-N}(F(y),\p(y'+\tau e_1-y))\neq0.
$$
Hence there exists a point $t'\in B(y'+\tau e_1,\a)\cap T$. For the same reason, there exists a point $t''\in B(y''+\tau e_1,\a)\cap T$ with the same $\tau$. We have
\begin{equation}\label{s2}
 |t'-t''|\le|t'-y'-\tau e_1|+|t''-y''-\tau e_1|+|y'-y''|<2\a+|y'-y''|<\eta(T).
\end{equation}
But this inequality is impossible.

Now suppose that our assumption is proved for all $n,\,m$ such that $n+m>M$ with $0\le M<2N$; then we prove it for $p, q\ge0$ such that $p+q=M$.

Assume the converse. There are points $y'\in\supp g_p,\ y''\in\supp g_q$ such that $0<|y'-y''|<\eta(T)$.  Pick
$$
\a<(1/2)\min\{|y'-y''|,\,(\eta(T)-|y'-y''|)\}
$$
and $C^\infty$-function $\p$ with $\supp\p\subset B(0,\a)$ such that
 $$
\P_{p,p}(y')=(g_p(y),\p(y'-y))\neq0,\qquad \P_{q,q}(y'')=(g_q(y),\p(y''-y))\neq0.
 $$
Take $\e>0$ such that $|\P_{p,p}(y')|>2\e$ and $|\P_{q,q}(y'')|>2\e$. Obviously, for a common $\e$-almost period $\tau$ of the functions
$\P_{p,p}(y'+se_1)$ and $\P_{q,q}(y''+se_1)$
  we get
  $$
  |\P_{p,p}(y'+\tau e_1)|>\e,\qquad |\P_{q,q}(y''+\tau e_1)|>\e.
  $$
   If there exists $n>p$ such that at least one of the numbers $\P_{n,j}(y'+\tau e_1),\,j=0,\dots,n$ does not vanish,
 then there is a point $t'\in B(y'+\tau e_1,\a)\cap\supp g_n$. Since $\P_{q,q}(y''+\tau e_1)\neq0$, we see that there is a point $t''\in B(y''+\tau e_1,\a)\cap\supp g_q$. By \eqref{s2}, $|t'-t''|<\eta(T)$, hence
 this inequality contradicts to our assumption $\eta(\supp g_n\cup\supp g_q)\ge\eta(T)$ for $n+q>M$. In the same way, we obtain a contradiction in the case
 when  at least one of the numbers $\P_{n,j}(y''+\tau e_1)$ for $n>q,\,j=0,\dots,n$,  does not vanish.

Therefore, for $t=y'+\tau e_1$ equality \eqref{s1} implies
$$
(y'_1+\tau)^{-p}(F(y),\p(y'+\tau e_1-y))=\P_{p,p}(y'+\tau e_1)+\sum_{m=0}^{p-1} (y'_1+\tau)^{m-p} \sum_{n=m}^{p-1}\P_{n,m}(y'+\tau e_1),
$$
and for $t=y''+\tau e_1$  implies
$$
(y''_1+\tau)^{-q}(F(y),\p(y''+\tau e_1-y))=\P_{q,q}(y''+\tau e_1)+\sum_{m=0}^{q-1} (y''_1+\tau )^{m-q} \sum_{n=m}^{q-1}\P_{n,m}(y''+\tau e_1).
$$
 If $\tau$ is large enough, then we get
$$
(y'_1+\tau)^{-p}(F(y),\p(y'+\tau e_1-y))\neq0,\qquad (y''_1+\tau)^{-q}(F(y),\p(y''+\tau e_1-y))\neq0,
$$
hence there exists $t'\in B(y'+\tau e_1,\a)\cap T$ and $t''\in B(y''+\tau e_1,\a)\cap T$ with the same $\tau$. By \eqref{s2}, we obtain a contradiction. Thus our assumption is valid for all $n,\,m\le N$. \bs

\medskip

{\bf Proof of Theorem \ref{T10}}. It follows from Proposition \ref{P1} i) that there is $K\in\N\cup\{0\}$ such that $\mu_k\equiv0$ for $\|k\|>K$.

In the case a), for any $\p\in S(\R^d)$ and $k$ we have
\begin{equation}\label{c}
(\hat\mu_k(y),\p(t-y))=(\mu_k(x),\check\p(x)e^{-2\pi i\la x,t\ra})=\sum_{\l\in\L}p_k(\l)\check\p(\l)e^{-2\pi i\la\l,t\ra}.
\end{equation}
Since $\check\p\in S(\R^d)$, we get $\check\p(x)=o(|x|^{-d-1})$. Hence,
$$
\sum_{\l\in\L}|p_k(\l)||\check\p(\l)|\le C_1+C_2\sum_{\l\in\L,|\l|>1}|\l|^{-d-1}=C_1+C_2\int_1^\infty s^{-d-1}dn(s),
$$
where $n(s)=\#(\L\cap B(0,s))$. The set $\L$ is of bounded density, therefore, $n(s)=O(s^d)$ and the integral
$$
 \int_1^\infty s^{-d-1}n(ds)=-n(1)+(d+1)\int_1^\infty\frac{n(s)ds}{s^{d+2}}
$$
is finite. Consequently the sum in \eqref{c} is absolutely convergent, and the function $(\hat\mu_k(y),\p(t-y))$ is almost periodic in $t\in\R^d$.

In the case b), by Proposition \ref{P1} ii), we get
\begin{equation}\label{p}
p_k(\l)=O(|\l|^M) \quad\mbox{as}\quad \l\to\infty.
\end{equation}
 Since $\check\p(x)=o(|x|^{-M-d-1})$ for each $H$, we can repeat the previous arguments and obtain
that $\hat\mu_k$ are almost periodic distributions.

Further, by \eqref{f2}, we have
$$
  \hat f=\sum_{k_1=0}^K (2\pi i)^{k_1}y_1^{k_1}g_{k_1}, \quad\mbox{where}\quad g_{k_1}(y)=\sum_{k_2+\dots+k_d\le K-k_1}(2\pi i)^{k_2+\dots+k_d}y_2^{k_2}\dots y_d^{k_d}\hat\mu_{k_1,\dots,k_d},
$$
and for any $\p\in\C^\infty$ with compact support
$$
(g_{k_1}(y),\p(t-y))=\sum_{k_2+\dots+k_d\le K-k_1}(2\pi i)^{k_2+\dots+k_d}(\hat\mu_{k_1,\dots,k_d},y_2^{k_2}\dots y_d^{k_d}\p(t-y)).
$$
Note that each term of the last sum can be rewritten as
$$
  \sum_{m_2\le k_2,\dots,m_d\le k_d}c_{m,k}t_2^{k_2-m_2}\dots t_d^{k_d-m_d}\left[(\hat\mu_k,\,(t_2-y_2)^{m_2}\dots (t_d-y_d)^{m_d}\p(t-y))\right]
$$
with some constants $c_{m,k}$. Since $\hat\mu_k$ are almost periodic distributions, we get that the expressions in square brackets are almost periodic functions in $t\in\R^d$, and hence in $t_1\in\R$.
Therefore the functions $(g_{k_1}(y),\p(t-y))$ are almost periodic in $t_1\in\R$ for any fixed $(t_2,\dots,t_d)\in\R^{d-1}$. Applying  Lemma \ref{L1} to distributions $g_{k_1},\,k_1=0,\dots,K$,
we get  they have uniformly discrete supports and $\eta(\supp g_{k_1})\ge\eta(\G)$.

 For a fixed $k_1$ we have
 $$
 g_{k_1}=\sum_{k_2=0}^{K-k_1}(2\pi i)^{k_2}y_2^{k_2} g_{k_1,k_2},
  $$
where
$$
 g_{k_1,k_2}(y)=\sum_{k_3+\dots+k_d\le K-k_1-k_2}(2\pi i)^{k_3+\dots+k_d}y_3^{k_3}\dots y_d^{k_d}\hat\mu_{k_1,\dots,k_d}.
$$
The functions $(g_{k_1,k_2}(y),\p(t-y))$ are almost periodic in $t_2\in\R$ for any fixed $(t_1,t_3\dots,t_d)\in\R^{d-1}$. Applying  Lemma \ref{L1} to distributions $g_{k_1,k_2},\,k_2=0,\dots,N-k_1$
with respect to the variable $y_2$, we get that these distributions have uniformly discrete supports and $\eta(\supp g_{k_1,k_2})\ge\eta(\G)$.

After a finite number of steps we get that  the support $\G_k$ of the distributions $\hat\mu_k$ for all $k\in(\N\cup\{0\})^d$ are
uniformly discrete and satisfy \eqref{et}. If the masses $p_k(\l)$ are uniformly bounded, then Proposition \ref{P1} iii), yields that $\hat\mu_k$ are measures with uniformly bounded masses.
If masses $p_k(\l)$ satisfy \eqref{p}, we have
$$
  \sum_{|\l|\le r}|p_k(\l)|\le C_1+C_2r^M n(r)=O(r^{M+d}) \quad\mbox{as}\quad r\to\infty.
$$
Applying Proposition \ref{P1} iii) with $K=0$,  we obtain that the coefficients $q_{k,j}(\g)$ in \eqref{m2} are uniformly bounded. \bs
\medskip

{\bf Remark 3}. Let $\L$ be uniformly discrete. If we replace the condition  $p_k(\l)=O(1)$ by $p_k(\l)=o(|\l|)$, then we get for any $\e>0$ and appropriate $r(\e)$
\begin{equation}\label{b1}
  \sum_{|\l|\le r}|p_k(\l)|\le \sum_{|\l|\le r(\e)}|p_k(\l)|+\sum_{r(\e)<|\l|\le r}|p_k(\l)| \le C(\e)+\e r\,n(r).
\end{equation}
Hence, \eqref{b1} is $o(r^{d+1})$ as $r\to\infty$. Applying Proposition \ref{P1} iii) with $K=0,\,M=1$ and taking into account Remark 1,
we obtain that $\hat\mu_k$ (and, of course, $\check\mu_k$) are measures with uniformly bounded masses. Applying
 Proposition \ref{P1} iii) to the measures $\check\mu_k$, we obtain that the coefficients $p_k(\l)$ are uniformly bounded.

\section{Proofs of the main theorems}\label{S7}

{\bf Proof of Theorem \ref{T8}}. By Theorem \ref{T10}, the measures $\mu_k$ from \eqref{c2} have uniformly discrete spectra. Applying Theorem \ref{T3} to these measures,
we get that $\supp\mu_k$ is contained in a finite union of translates of a full-rank lattice.
 To obtain the last assertion of the theorem,  we  take the measures $e^{-i\xi}\mu_k$ instead of the measures $\mu_k$.  \bs
\medskip

{\bf Proof of Theorem \ref{T9}}.
Notice that   masses of the measures $\mu_k$ from \eqref{c2} are not assumed to be separated from zero, therefore we cannot make use of the result of Theorem \ref{T5}.

Let $\p(y)$ be nonnegative $C^\infty$-function with support in the unit ball such that \newline $\int\p(t)dt=1$. Pick $k\in(\N\cup\{0\})^d$, and set $\rho_\e=\e^d\p(\e y)\check\mu_k(y)$,
where $\check\mu_k$ is the inverse Fourier transform of the measure $\mu_k$. We have
\begin{equation}\label{j}
 \hat\rho_\e(s)=\int_{\R^d}\hat\p\left(\frac{s-x}{\e}\right)\mu_k(dx)=\left(\int_{|s-x|<1}+\int_{|s-x|\ge1}\right)\hat\p\left(\frac{s-x}{\e}\right)\mu_k(dx).
\end{equation}
Since $\hat\p\in S(\R^d)$, we get $\hat\p(y)=o(|y|^{-d-1})$ as $|y|\to\infty$ and $\hat\p((s-x)/\e)\to0$ for $s\neq x$ as $\e\to0$. Also, $\hat\p(0)=1$. Hence, by the Dominated Convergence Theorem,
$$
\lim_{\e\to0}\int_{|s-x|<1}\hat\p\left(\frac{s-x}{\e}\right)\mu_k(dx)=\begin{cases}\mu_k(\{\l\})=p_k(\l), & s=\l\in\L,\\  0, & s\not\in\L.\end{cases}
$$
Note that the masses $\mu_k(\{\l\})$ are uniformly bounded. Also,
$$
n_s(t):=\#(\L\cap B(s,t))=O(t^d)\qquad(t\to\infty).
$$
 Therefore,
$$
\left|\int\limits_{|s-x|\ge1}\hat\p\left(\frac{s-x}{\e}\right)\mu_k(dx)\right|\le C\int\limits_1^\infty(t/\e)^{-d-1}n_s(dt)\le\e^{d+1}(d+1)C\int\limits_1^\infty t^{-d-2}n_s(t)dt,
$$
and the left-hand side tends to zero as $\e\to0$.

Furthermore, by Theorem \ref{T10}, $\hat\mu_k$ has uniformly bounded masses and support of bounded density. Hence, $|\hat\mu_k|(B(0,r))=O(r^d)$ as $r\to\infty$ and the same is valid for the measures $\check\mu_k$.
Therefore,
$$
|\rho_\e|(B(0,r))\le\e^d|\check\mu_k|(B(0,1/\e)=O(1)\qquad (\e\to0).
$$
Since support of the measure $\rho_\e$ is bounded for each fixed $\e$, we see that the measures $\rho_\e$ form a bounded family of linear functionals on the space of continuous bounded functions on $\R^d$
 and, in particular, on its subspace $C(\hR)$. A closed  ball in the conjugate space $C^*(\hR)$ is a compact set in the weak-star topology. Hence there exists a measure $\hr_k$ on $\hR$ with a total variation
 $\|\hr_k\|<\infty$ such that for every $f\in C(\hR)$ there is a subsequence $\e'\to0$, for which we have $\la\rho_{\e'},f\ra\to\la\hr_k,f\ra$.
 Applying this to every character $x\in\R^d$ in the place of $f$, we obtain from \eqref{j}
 $$
 \hat\hr_k(x)=\lim_{\e'\to0}\hat\rho_{\e'}(x)=\begin{cases}p_k(\l), & x=\l\in\L,\\  0, & x\not\in\L.\end{cases}
 $$
 Note that $\hat\hr_k(x)$ is a continuous function with respect to the discrete topology on $\R^d$, and $|p_k(\l)|\le\|\hr_k\|$ for all $\l\in\R^d$.

Let $\psi$ be a $C^\infty$-function with compact support.  Since the inverse Fourier transform of the function $\psi(t-x)$ is $e^{2\pi i\la t,y\ra}\hat\psi(y)$, we get
\begin{equation}\label{f0}
(\psi\star\mu_k)(t)=\int_{\R^d}\psi(t-x)\mu_k(dx)=\int_{\R^d}e^{2\pi i\la y,t\ra}\hat\psi(y)\hat\mu_k(dy)=\sum_{\g\in\supp\hat\mu_k}\hat\mu_k(\g)\hat\psi(\g)e^{2\pi i\la\g,t\ra}.
\end{equation}
Clearly, $\hat\psi\in S(\R^d)$ and $\hat\psi(\g)=o(|\g|^{-d-1})$ as $\g\to\infty$. Also, the measure $\hat\mu_k$ has uniformly bounded masses and uniformly discrete support, therefore the series in the right-hand side is absolutely converges and $(\psi\star\mu_k)(t)\in W$. In particular, the measure $\mu_k$ is almost periodic.

 Furthermore, the function
$$
F(t)=\sum_{\|k\|\le K}(\psi\star\mu_k)(t)\overline{(\psi\star\mu_k)(t)}
$$
belongs to $W$ as well. Now suppose that $\psi$ has support in the ball $B(0,\eta(\L)/2)$ and $\psi(0)=1$. Taking into account that
\begin{equation}\label{p0}
(\psi\star\mu_k)(\l)=p_k(\l)\qquad\mbox{for}\quad \l\in\L,
\end{equation}
and using \eqref{c3}, we get
$$
 F(\l)=\sum_{\|k\|\le K} |p_k(\l)|^2\ge\a>0\qquad \forall\, \l\in\L.
$$
Using Proposition \ref{P3} with $h(z)=1/z$, we  construct the function
$$
g(t)=\sum_nc_ne^{2\pi i\la t,\tau_n\ra}\in W
$$
 with the property $g(\l)=1/F(\l)$ for $\l\in\L$.

Consider the dual pair $(\hR,\R^d_{discr})$. Set
$$
\ha=\sum_{\|k\|\le K}\hr_k(\t)\star\overline{\hr_k(-\t)},\quad \t\in\hR.
$$
 We get
 $$
 \hat\ha(\l)=F(\l)\quad \forall \l\in\L,\quad\mbox{and}\quad  \hat\ha(x)=0\quad\forall x\not\in\L.
 $$
Since  $\sum_n|c_n|<\infty$, we see that the series $\sum_nc_n\ha(y+\tau_n)$ converges with respect to the norm to a measure $\hb(y)$ such that
 $$
 \hat\hb(x)=\sum_n c_ne^{2\pi i\la x,\tau_n\ra}\hat\ha(x)=\begin{cases}1, & x\in\L,\\ 0, & \mbox{otherwise}. \end{cases}
 $$
By Theorem \ref{T6}, $\L$ belongs to the coset ring of $\R^d_{discr}$. By Theorem \ref{T7},
\begin{equation}\label{b2}
\L=\cup_{n=1}^N(A_n\setminus\cup_{s=1}^{S_n} B^s_n),
\end{equation}
where $A_n,\,B^s_n$ are cosets of discrete lattices.

Let $A,\,A'$  be two arbitrary cosets in $\R^d$ and $L,\,L'$ be the corresponding lattices. Clearly, $A\cap A'$ is a coset of $L\cap L'$. If $r(A\cap A')=d$, then the factor-groups $L/(L\cap L')$ and $L'/(L\cap L')$ are finite.
Therefore, $A\setminus A'=A\setminus(A\cap A')$ is a finite union of disjoint cosets of $L\cap L'$, and the same is valid for
$$
A\cup A'=(A\setminus(A\cap A'))\cup(A'\setminus(A\cap A'))\cup(A\cap A').
$$
We repeat this transformation for each pair of cosets with dimension $d$ of their intersection. After a finite number of steps, we get a representation similar to \eqref{b2} such that  $r(B^s_n)<d$
and $r(A_n\cap A_{n'})<d$ for $n\neq n'$.

Put $\d_A=\sum_{x\in A}\d_A$. Using the equalities $\d_{A\setminus A'}=\d_A-\d_{A\cap A'}$ and $\d_{A\cup A'}=\d_A+\d_{A'}-\d_{A\cap A'}$, we get
$$
   \d_\L=\sum_{j=1}^J\d_{D_j}+\sum_{j=1}^{J_1}\d_{E_j}-\sum_{j=1}^{J_2}\d_{F_j},
$$
where  $D_j,\,E_j,\,F_j$ are cosets such that $r(D_j)=d$ and $r(E_j)<d,\,r(F_j)<d$.

 Every coset $D_j$ is a translation of full-rank lattice,  hence $\d_{D_j}$ is $d$-periodic and $\sum_{j=1}^J\d_{D_j}$ is almost periodic.
Above we showed that all measures $\mu_k$ are almost periodic. By Proposition \ref{P2} ii), the measures $|\mu_k|$ are almost periodic, hence the measure $\sum_k|\mu_k|$
with support $\L$ is almost periodic as well.
 Using \eqref{c3} and Proposition \ref{P2} iii), we get that the measure $\d_\L$ is almost periodic. Therefore the left-hand side of the equality
 \begin{equation}\label{b}
 \sum_{j=1}^{J_1}\d_{E_j}- \sum_{j=1}^{J_2}\d_{F_j}=\d_\L-\sum_{j=1}^J\d_{D_j}
 \end{equation}
is an almost periodic measure as well. But its support  is contained in a finite union of hyperplanes and isn't relatively dense. Therefore the measure in the left-hand side of \eqref{b} is identically zero, and $\d_\L=\sum_{j=1}^J\d_{D_j}$.
Since  the measure $\d_\L$ has the unit masses at every point of $\L$, we see that the cosets $D_j$ are pairwise disjoint and  for some full-rank lattices $L_j$ and $\l_j\in\R^d$
\begin{equation}\label{m3}
\L=\cup_{j=1}^J D_j=\cup_{j=1}^J(\l_j+L_j).
\end{equation}

Recall that $\psi\star\mu_k\in W$. It follows from \eqref{p0} and \eqref{m3} that
\begin{equation}\label{m0}
 \mu_k=\sum_{j=1}^J (\psi\star\mu_k)\d_{L_j+\l_j}=\sum_{j=1}^J\sum_s a_{s,k}e^{2\pi i\la x,\g_{s,k}\ra}\d_{L_j+\l_j}
\end{equation}
for some $\g_{s,k}\in\R^d$ and $a_{s,k}\in\C$ such that $\sum_s|a_{s,k}|<\infty$.
  For every $\g\in\R^d$ there is $\b$ inside the parallelepiped generated by corresponding $L_j^*$ such that $\g-\b\in L_j^*$. Therefore,
$e^{2\pi i\la x,\g_{s,k}\ra}=e^{2\pi i\la x,\b_{s,k,j}\ra}$ for $x\in L_j$, and the set of all points $\b_{s,k,j}$ is bounded.
Therefore \eqref{m0} can be rewritten as
$$
\sum_{j=1}^J\sum_s\sum_{x\in L_j} b_s(k,j)e^{2\pi i\la x,\b_{s,k,j}\ra}\d_{x+\l_j}\quad\mbox{with}\quad b_s(k,j)= a_{s,k}e^{2\pi i\la \l_j,\g_{s,k}-\b_{s,k,j}\ra}.
$$
By \eqref{f2}, we obtain the statement of Theorem \ref{T9}.  \bs

\medskip

Here we actually proved the following statement:
\begin{Pro}\label{P4}
Let $\{\nu_n\}$ be a finite set of measures from $S^*(\R^d)$ such that for each $n$

a) $\supp\nu_n$ is uniformly discrete,

b) $\supp\hat\nu_n$ is locally finite of bounded density,

c) $\hat\nu_n$ is a measure with uniformly bounded masses.
\smallskip

\noindent Set $\hL=\cup_n\supp\nu_n$. If
$$
\inf_{\l\in\hL}\max_n|\nu_n(\l)|>0,
$$
then $\hL$ is  a  finite union of translates of several full-rank lattices.
\end{Pro}
Note that the boundedness of masses $\nu_n(\l)$ follows from \eqref{f0} with $\nu_n$ instead of $\mu_k$ and $\psi$ such that $\supp\psi\subset B(0,\eta(\supp\nu_n))$ and $\psi(0)=1$.

{\bf Remark 4}. It is not hard to check that conditions b) and c) can be replaced by

b') $\hat\nu_n$ is a purely point measure,

c')  $|\hat\nu_n|(B(0,r))=O(r^d)$ as $r\to\infty$.

\medskip

\bigskip
I am very grateful to the referee for carefully reading my article and giving numerous remarks.


\begin{thebibliography}{8}

\bibitem{A} L.N.Argabright and J.G.de Lamadrid, {\it Almost Periodic Measures}, Memoirs of the
American Mathematical Society {\bf 428} (1990).

\bibitem{D} M.Baake, R.Moody (eds.), {\it Directions in Mathematical Quasicrystals}, 
CRM Monograph Series, Vol. 13, American Mathematical Society, Providence, RI, 2000.

\bibitem{Co}  A.Cordoba,  {\it Dirac combs}, Letters in Mathematical Phisics {\bf 17} (1989), 191-196.

\bibitem{C}  C.Corduneanu, {\it  Almost Periodic Functions},  Chelsea, New-York, 1989. 

\bibitem{F2} S.Yu.Favorov, {\it Fourier Quasicrystals and Lagarias' Conjecture},
 Proceedings of the American Mathematical Society {\bf 144} (2016) , 3527-3536.

\bibitem{F3} S.Yu.Favorov, {\it Some Properties of Measures with Discrete Support}, Matematychni Studii  {\bf 46} (2016), 189-195.

\bibitem{F1} S.Yu.Favorov, {\it Tempered distributions with discrete support and spectrum}, Bulletin of the Hellenic Mathematical Society {\bf 62} (2018), 66-79.

 \bibitem{F4} S.Yu.Favorov, {\it Large Fourier Quasicrystals and Wiener's Theorem},  Journal of Fourier Analysis and Applications {\bf 25} (2019), 377-392.

\bibitem{F5} S.Yu.Favorov, {\it Local Wiener's Theorem and Coherent Sets of Frequencies},  Analysis Mathematica {\bf 46} (2020), 737–746.

\bibitem{FK} S.Yu.Favorov  and Ye.Kolbasina, {\it Perturbations of discrete lattices and almost
periodic sets}, Algebra and Discrete Mathematica {\bf 9} (2010), 48-58.

\bibitem{K2} M.N.Kolountzakis, {\it Fourier Pairs of Discrete Support with Little Structure},
Journal of Fourier Analysis and Applications {\bf 22} (2016), 1-5.

\bibitem{K1} M.N.Kolountzakis, {\it On the Structure of Multiple Translations Tilings by Polygonal Regions},
https://arxiv.org/abs/math/9904065.

\bibitem{KL} M.N.Kolountzakis and J.C.Lagarias, {\it Structure of Tilings of the Line by a Function},
Duke Math.Journal {\bf 82} (1996), 653-678.

\bibitem{KS} P.Kurasov and R.Suhr, {\it Asymptotically isospectral quantum graphs and generalised trigonometric polynomials}, 
Journal of Mathematical Analysis and Applications {\bf 488} (2020), Article no. 124049.

\bibitem{L2} J.C.Lagarias, {\it Geometric Models for Quasicrystals I.Delone Set of Finite Type},
Discrete \& Computational Geometry {\bf 21} (1999),161-191.

\bibitem{L1} J.C.Lagarias {\it Mathematical Quasicrystals and the Problem of Diffraction}, CRM Monograph Series, Vol. 13, American Mathematical
Society, Providence, RI, 2000, pp. 61–93.

\bibitem{LM} B.Leptin and D.M$\ddot{u}$ller, {\it Uniform partitions of unity on locally compact groups},
 Advances in Mathematics {\bf 90} (1991), 1-14.

\bibitem{LO1} N.Lev and A.Olevskii, {\it Measures with Uniformly Discrete Support and Spectrum},Comptes Rendus Math´ematique. Acad´emie des Sciences. 
Paris {\bf 351} (2013), 599–603.

\bibitem{LO3} N.Lev and A.Olevskii,  {\it Quasicrystals with Discrete Support and Spectrum},
 Revista Matem´atica Iberoamericana {\bf 32}(2016), 1341-1252.

\bibitem{LO2} N.Lev and A.Olevskii, {\it Fourier Quasicrystals and Discreteness of the Diffraction Spectrum}, Advances in Mathematics {\bf 315} (2017), 1-26.

\bibitem{LR}  N.Lev and G.Reti, {\it Crystalline Temperate Distribution with Uniformly Discrete Support and Spectrum}, 
Journal of Functional Analysis {\bf 281} (2021), Article no. 109072.

\bibitem{M3} Y.Meyer, {\it Nombres de Pisot, Nombres de Salem et analyse harmonique}, Lecture Notes
in Mathematics, Vol. 117, Springer, Berlin–New York, 1970.

\bibitem{M1} Y.Meyer,  {\it Quasicrystals, Almost Periodic Patterns, Mean--periodic Functions, and Irregular Sampling}, 
African Diaspora Journal of Mathematics {\bf 13} (2012), 1-45.

\bibitem{M2} Y.Meyer, {\it Guinand's Measure are  Almost Periodic Distributions},
Bulletin of the Hellenic Mathematical Society {\bf 61} (2017), 11-20.

\bibitem{P} V.P.Palamodov, {\it A Geometric Characterization of a Class of Poisson Type Distributions}, Journal of Fourier Analysis and Applications
{\bf 23} (2017), 1227-1237.

\bibitem{Q} J.Patera (ed.), {\it Quasicrystals and Discrete Geometry}, Fields Institute Monographs,
Vol. 10, American Mathematical Society, Providence RI, 1998.

\bibitem{R} L.I.Ronkin, {\it Almost Periodic Distributions and Divisors in Tube
Domains}, Zapiski Nauchnykh Seminarov POMI {\bf 247} (1997), 210–236.

\bibitem{Ru} W.Rudin, {\it Functional Analysis}, McGraw-Hill Series in Higher Mathematics, McGraw-Hill, New York–D¨usseldorf–Johannesburg, 1973.

\end{thebibliography}
\end{document}